\documentclass[12pt]{article}
\usepackage{graphicx}
\usepackage{amsfonts}
\usepackage{bbm}
\usepackage{mathrsfs}
 \usepackage{mathrsfs,amsmath,amssymb,amsthm}
\usepackage[dvips]{color}
 \usepackage{amssymb}
 \usepackage{amsbsy}
 \usepackage[dvips]{color}
\allowdisplaybreaks
 \theoremstyle{plain}
\theoremstyle{remark}  \newtheorem{remark}{\noindent\mbox{Remark}}
 \theoremstyle{plain}
 \theoremstyle{plain}\newtheorem{lemma}{\noindent\mbox{Lemma}}
\theoremstyle{plain} \newtheorem{theorem}{\noindent\mbox{Theorem}}
 \theoremstyle{plain}\newtheorem{proposition}{\noindent\mbox{Proposition}}
 \theoremstyle{plain}\newtheorem{corollary}{\noindent\mbox{Corollary}}
\theoremstyle{definition}

 \def\bq{\begin{equation}}
 \def\eq{\end{equation}}
 \def\eqn{\end{eqnarray}}
 \def\bqn{\begin{eqnarray}}
 \def\proof{\noindent{\it Proof.~~}}
 \def\qed{\hfill$\Box$\medskip}
 \def\rto{\rightarrow\infty}
 \def\z{\left}
 \def\y{\right}
 \def\no{\nonumber}
\def\mbz{\mathbb{Z}}
 \def\mbe{\mathbb{E}}
 \def\mbp{\mathbb{P}}
 \def\ka{{\kappa}}
 \def\mcp{\mathcal{P}}

\begin{document}
\begin{center}{\Large Law of large numbers for random walk with unbounded jumps and BDP with bounded jumps in random environment}\footnote{Supported by National
Nature Science Foundation of China (Grant No. 11226199) and Nature Science Foundation of Anhui Educational Committee (Grant No. KJ2014A085)}\end{center}

\vspace{0.5cm}
\centerline {Hua-Ming \uppercase{Wang}$^a$ }

    \begin{center}
      {\footnotesize$^a$Department of Mathematics, Anhui Normal University, Wuhu 241003, China

    E-mail\,$:$ hmking@mail.ahnu.edu.cn}
    \end{center}
%\date{}
%\maketitle%

\vspace{-.3cm}

\begin{center}
\begin{minipage}[c]{12cm}
\begin{center}\textbf{Abstract}\quad \end{center}

We study random walk with unbounded jumps in random environment. The environment is stationary and ergodic, uniformly elliptic and decays polynomially with speed $Dj^{-(3+\varepsilon_0)}$ for some small $\varepsilon_0>0$ and proper $D>0.$ We prove a law of large number with positive velocity under the condition that the annealed mean of the hitting time of the positive half lattice is finite. Secondly, we consider birth and death process with bounded jumps in stationary and ergodic environment. Under the uniformly elliptic condition, we prove a law of large number and give the explicit formula of its velocity.

\vspace{0.2cm}


\vspace{0.2cm}
\end{minipage}
\end{center}

\section{Introduction}

Our prime concern is to prove the Law of Large Number (LLN hereafter) for Birth and Death Process in Random Environment (BDPRE hereafter) with bounded jumps , say $\{N_t\}_{t\ge0}.$ We assume that at each discontinuity, the particle jumps at most a distance $R$ to the right or at most a distance $L$ to the left. For the nearest neighbour setting $L=R=1,$ if one defines $T_n$ the hitting time of $n,$ then $N_{T_n}=n$ and $\{T_n-T_{n-1}\}_{n\ge1}$ forms a stationary and mixing sequence under the annealed probability whenever the environment is i.i.d.. Therefore the LLN of $\{N_t\}_{t\ge0}$ follows from that of $\{T_n\}_{n\ge0}.$ For details, see Ritter \cite{rit}. For bounded-jump setting, the above approach does not work. Similarly, define the ladder times $T_0=0$ and $T_n=\inf\{t>T_{n-1}:N_{T_n}>N_{T_{n-1}}\},\ n\ge1.$ Two problems arise: (1) It is hard to tell the exact value of $N_{T_n},$ though we know that $n\le N_{T_n}\le nR;$ (2) $\{T_n-T_{n-1}\}_{n\ge1}$ is not a stationary sequence under the annealed probability even if the environment is i.i.d..

We turn to consider the $h$-skeleton process $\{N_{nh}\}_{n\ge0},\ h>0$ of $\{N_t\}_{t\ge0}.$ If the LLN of the $h$-skeleton process is proved, then the LLN of $\{N_t\}_{t\ge0}$ follows from some standard procedures. However, $\{N_{nh}\}_{n\ge0}$ is indeed a discrete time {\it random walk in random environment with unbounded  jumps,} since theoretically speaking, $\{N_t\}_{t\ge0}$ may have many jumps in a time interval of length $h.$

Although the developments of Random Walk in Random Environment (RWRE hereafter) with bounded jumps  were almost satisfying, RWRE with unbounded jumps was very seldom considered.  Andjel \cite{adj} proved the 0-1 law. Comets and Popov \cite{cp} proved  an LLN. In \cite{cp}, two main conditions were required essentially: a) the jumping probabilities of the walk have an exponential tail; b)   `` the random
walk is `uniformly' transient to the right (i.e., there are no `traps')".
We mention also that in Gallesco and Popov \cite{gpa,gpb}, the author studied the central limit theorem of random walks with unbounded jumps among random conductances.

The BDPRE with bounded jumps considered in this paper is the continuous time analogue of RWRE with bounded jumps.
The situation of RWRE with bounded jumps on $\mathbb Z$ is almost satisfying. We review here only the known results closely related to the transient RWRE. For the nearest neighbour setting, the 0-1 law and the LLN with explicit velocity in the case of i.i.d. environment were studied in Solomon \cite{sol}; the extension to stationary and ergodic environment
was given in Alili \cite{ali}. The RWRE with bounded jumps (jumps are non-nearest neighbor) was introduced in Key \cite{key}, where the middle Lyapunov exponents of a sequence of positive random matrices were used to give the recurrence criteria.
Letchikov \cite{letc} simplified the recurrence criteria and proved a $\log^2n$ law for the recurrent case as Sinai's walk.
The
regime where the LLN holds with a positive velocity was characterized
in Br\'emont \cite{brb,brc}  by ``the environment viewed from particle",
but no explicit formula is  available any longer. In \cite{hwa,hwb,hz}, the authors set up the branching structure for RWRE with bounded jumps and gave  the explicit formulae for the velocities of LLNs. At last we point out that Bolthausen and Goldsheid \cite{bg} studied the recurrence and transience of RWRE on a strip, which is a generalization of RWRE with bounded jumps.
The LLN for RWRE on a strip could be found in Roitershtein \cite{roi}.

In this paper, firstly we  prove an LLN for RWRE with unbounded jumps. In our setup, we only need to require that the jumping probabilities have an polynomial tail. Moreover  we do not need the condition `` the random walk is `uniformly' transient to the right" used in \cite{cp}. Our approach is based on ``the environment viewed from particle" which dates back to Kozlov \cite{K85}. But in Kozlov \cite{K85} and some later literatures, this approach is only used to treat the bounded-jump setting. We use the large deviation for martingales (with unbounded martingale differences) derived in Lesigne and Voln\'y \cite{lv} to estimate the tail of a martingale related to RWRE with unbounded jumps. By constructing an invariant measure for the auxiliary Markov chain of the environment viewed from particle, we could show the LLN for the summation of the local drift of RWRE with unbounded jumps.  In this way, we  prove the LLN for  RWRE with unbounded jumps.
Secondly, using the LLN of RWRE with unbounded jumps, we prove the LLN of BDPRE with bounded jumps. We assume a uniform ellipticity  condition on the environment, which is crucial for our development. Under such  ellipticity condition, we could show that the tail probability of the jumps for the $h$-skeleton process is exponentially bounded. Therefore we could use the LLN for RWRE with unbounded jumps to get the LLN for the $h$-skeleton process. Then the LLN for BDPRE with bounded jumps follows by some standard procedures.

The paper is organized as follows. In Section \ref{mr}, we define strictly the models  and state the main results. The LLN for RWRE with unbounded jumps is proved in Section \ref{pr}  while
the  LLNs for the  $h$-skeleton process $\{N_{nh}\}_{n\ge 0}$ and $\{N_t\}$ are proved in Section \ref{pec}.
We devote Section \ref{dvp} to give the explicit formula of the velocity $v_{\tilde\mbp}$ for the case $L=R=2$ by using the branching structure in the embedded process constructed in \cite{hwb}. An appendix section is also given at the end of the paper to discuss the existence of the process $\{N_t\}.$

\section{The models and results}\label{mr}

\subsection{RWRE with unbounded jumps}
Let $\Omega$ be the collection of $\omega=(\omega_{x})_{x\in\mathbb Z}$ where for $x\in \mathbb Z,$ $\omega_x=(\omega_{xy})_{y\in\mathbb Z}$ is a probability measure on $\mathbb Z,$ that is, $\omega_{xy}\ge 0$ for all $y\in\mathbb Z$ and $\sum_{y\in\mathbb Z}\omega_{xy}=1.$   Let $\theta$ be the shift operator on $\Omega$ defined by $(\theta\omega)_x:=\omega_{x+1}.$ Equip $\Omega$ with Borel $\sigma$-algebra $\mathcal F$ and let $\mbp$ be a probability measure on $(\Omega,\mathcal F).$ For a typical realization of $\omega,$ we consider a Markov chain $\{S_n\}_{n\ge 0}$ with transitional probabilities $$P_{\omega}^{x_0}(S_{n+1}=x+y\big |S_{n}=x)=\omega_{xy}\text{ for all }n\ge0,\ P_\omega^{x_0}(S_0=x_0)=1,$$
so that $P_\omega^{x_0}$ is the quenched law of the Markov chain starting from $x_0$ in the
environment $\omega.$
Define a new probability measure $P^{x_0}$ by $$P^{x_0}(\cdot)=\int P^{x_0}_\omega(\cdot)\mbp(d\omega),$$ which is called the annealed probability. We use $E_\omega^{x_0},$ $E^{x_0}$ and $\mbe$ to denote the expectation operator for $P_\omega^{x_0},$ $P^{x_0}$ and $\mbp$ respectively. The superscript $x_0$ will be omitted whenever it is $0.$

\noindent{\bf Condition B}

\noindent\textbf{(B1)} $(\Omega,\mathcal F, \mathbb P, \theta)$ forms a stationary and ergodic system.

\noindent\textbf{(B2)} There exists $\varepsilon>0$ such that $\mbp(\omega_{01}>\varepsilon)=1.$

\noindent\textbf{(B3)} There exist small $\varepsilon_0>0$ and proper $D>0,$ such that $\mbp$-a.s., $$\omega_{0j}<D|j|^{-(3+\varepsilon_0)}.$$

Define $T=\inf\{n>0:S_n>0\},$ which is the time $\{S_n\}_{n\ge 0}$ hits $[1,\infty)$  and let  $$U_k=\#\{0\le n<T:S_n=k\}.$$
Here and throughout, $``\#\{\ \}"$ denotes the number of elements  in set $\{\ \}.$
  We have the following ballistic LLN for $\{S_n\}_{n\ge0}.$
\begin{theorem}\label{unb}
  Suppose that condition B holds and $E(T)<\infty.$ Then  $$P\text{-a.s.,}\ \lim_{n\rto}\frac{S_n}{n}=v_\mbp>0$$
  where
   $$
  v_\mbp=\frac{\mbe\z(\sum_{i=1}^\infty\sum_{k\le 0}E_{\theta^{-k}\omega}\big(U_k|S_{T=i}\big)\sum_{j\in\mbz} j\omega_{0j}\y)}{\sum_{i=1}^\infty E\big(T|S_{T}=i\big)}.
$$
\end{theorem}
\begin{remark}
(1) Our proof is based on an approach known as ``the environment viewed from particle" introduced in Kozlov \cite{K85}. In \cite{K85} and the later literatures, this approach was only used to  treat the bounded-jump setting. In our setup, we construct an invariant measure $Q$ which is equivalent to $\mbp.$  Under $Q,$ $\overline{\omega}(n):=\theta^{S_n}\omega,\ n\ge0$ form a stationary and ergodic sequence.
Let $d(x,\omega)=E_\omega^x(S_1-S_0)$  be the local drift and set $$M_n=S_n-S_0-\sum_{k=0}^{n-1}d(S_n,\omega).$$ $\{M_n\}_{n\ge0}$ is a martingale with unbounded differences. But under (B3), we could use the large deviation for martingale to estimate the tail of $M_n.$ Therefore, using the auxiliary Markov chain $\{\overline{\omega}(n)\}_{n\ge0}$ we prove the LLN of $\{S_n\}_{n\ge0}.$

(2) For RWRE with bounded jumps. The velocity of the LLN could be written in terms of the environment $\omega$ explicitly because one could calculate the quenched mean $E_\omega (T)$ by using the branching structure constructed from the path of the random walk. See \cite{hwb} for details. However for unbounded-jump setting, we do not know how to calculate the quenched mean $E_\omega( T).$ Therefore, the condition $``E(T)<\infty"$ looks not so satisfying and consequently the velocity could not be given explicitly. For some special case, we could give the explicit formula for the velocity $v_\mbp.$ See the discussion below.
\end{remark}

If $\mbp(\omega_{0j}=0,j\ge 2)=1,$  define $$\Phi=\left(
        \begin{array}{ccccccc}
          0 & 1 & 0 & 0 &0 & 0&\cdots \\
          \omega_{0,1} & \omega_{0,0} &  \omega_{0,-1}& \omega_{0,-2} &\omega_{0,-3}  &\omega_{0,-4}& \cdots \\
          0 & \omega_{-1,1} & \omega_{-1,0} & \omega_{-1,-1} & \omega_{-1,-2} & \omega_{-1,-3}&\cdots \\
          0 & 0 & \omega_{-2,1} & \omega_{-2,0} & \omega_{-2,-1} &\omega_{-2,-2}&\cdots \\
          0 & 0 & 0& \omega_{-3,1} & \omega_{-3,0} &\omega_{-3,-1} &\cdots \\
          0 & 0 & 0&  0&\omega_{-4,1}&\omega_{-4,0}&\cdots\\
          \vdots & \vdots & \vdots & \vdots & \vdots & \vdots&\ddots \\
        \end{array}
      \right).
$$
Let $\pi=(\pi_1,\pi_0,\pi_{-1},\pi_{-2},...)$ be a solution of the equation $$\pi\Phi=\pi,$$ with $\pi_1=1.$
If $\mbe(\sum_{i\le0}\pi_i)<\infty,$ then it follows from the classical ergodic theory for Markov chain that $P$-a.s., $T<\infty$ and $E_\omega(T)=\sum_{i\le 0}\pi_i.$ We have the following corollary of Theorem \ref{unb}.

\begin{corollary}
  Suppose that condition B holds and $\mbe (\sum_{i\le 0}\pi_i)<\infty.$ Then $$P\text{-a.s.,}\ \lim_{n\rto}\frac{S_n}{n}=\frac{1}{\mbe\z(\sum_{i\le0}\pi_i\y)}.$$
\end{corollary}

\subsection{BDPRE with bounded jumps}
Next we define the birth and death process in random environment with bounded jumps. To construct the environment, fix $1\le L,R\in \mathbb Z$ and let $\tilde\Omega$ be the collection of $\tilde\omega=(\tilde\omega_i)_{i\in\mathbb Z}=(\mu^{L}_i,...,\mu^{1}_i,\lambda^1_i,...,\lambda^R_i)_{i\in \mathbb Z},$ where $\mu_i^l,\lambda_i^r\ge 0$ for all $i\in\mathbb Z,\ l=1,..,L$ and $r=1,...,R. $ Equip $\tilde\Omega$ with the Borel $\sigma$-algebra $\tilde{\mathcal{F}}$ and let $\tilde{\mathbb P}$ be a probability measure on $(\tilde\Omega,\tilde{\mathcal F}).$ Then the so-called  random environment $\tilde\omega$ is a random element of $\tilde\Omega$ chosen according to $\tilde{\mathbb P}.$
Given a realization of $\tilde\omega,$ let $\{N_t\}_{t\ge0}$ be a continuous time Markov chain, which waits at a state $n$ an exponentially distributed time with parameter $\sum_{l=1}^L\mu^{l}_n+\sum_{r=1}^R\lambda^{r}_n$ and then jumps to $n-i$ with probability ${\mu^i_n}/(\sum_{l=1}^L\mu^{l}_n+\sum_{r=1}^R\lambda^{r}_n),$ $i=1,...,L$ or to $n+j$ with probability ${\lambda^j_n}/(\sum_{l=1}^L\mu^{l}_n+\sum_{r=1}^R\lambda^{r}_n),$ $j=1,...,R.$ In this paper, we always assume that the paths of $\{N_t\}$ are right continuous. We call the process $\{N_t\}_{t\ge 0}$ a {\it birth and death process in random environment with bounded jumps}.

%
%
%To figure out the explicit asymptotic velocity of LLN, we use the branching structure within the embedded process revealed in Hong and Wang \cite{hwb}. The idea is as follows. The invariant density for the ``environment viewed from particle" is closely related to the quenched mean of the ladder time $T_1:=\inf[t>0:N_t>0].$ By the branching structure within the embedded process $\{\chi_n\},$ one could use a multitype branching process in random environment to count how many times $\{N_t\}$ has ever visited state $i$ before $T_1.$ But every time it visits $i,$ it would wait here an exponentially distributed time period. In this way, we could study the distribution of $T_1$ and consequently give the explicit velocity for LLN of $\{N_t\}$ when $L=R=2.$

For a typical realization of $\tilde\omega,$ $\tilde P^x_{\tilde\omega}$ denotes the law induced by the process $\{N_t\}$ starting from $x.$ The measure $\tilde P^x_{\tilde\omega}$ is usually related as the {\it quenched} probability. Define the annealed probability measure $\tilde P^x$ by $\tilde P^x(\cdot)=\int_{\Omega}\tilde P^x_{\tilde\omega}(\cdot)\tilde{\mathbb P}(d\tilde\omega).$
The notations $\tilde E^x_{\tilde\omega},$ $\tilde E^x$ and $\tilde{\mathbb E}$ will be used to denote the expectation operators for $\tilde P^x_{\tilde\omega},$ $\tilde P^x$ and $\tilde{\mathbb P}$ respectively. The superscript $x$  will be omitted if it is $0.$ Let operator $\theta$ be the canonical shift on $\tilde\Omega$ defined by $(\theta\tilde\omega)_i=\tilde\omega_{i+1}.$

\noindent {\bf Condition C}

\noindent\textbf{(C1)} $(\tilde\Omega,\tilde{\mathcal F}, \tilde{\mathbb P}, \theta)$ forms a stationary and ergodic system.

\noindent\textbf{(C2)} the measure $\tilde{\mathbb P}$ is uniformly elliptic, that is, $$\tilde{\mathbb P}\Big(\varepsilon<\mu_0^l, \lambda_0^r<M , 1\le l\le L, 1\le r\le R \Big)=1$$ for some small $\varepsilon>0$ and large $M>0.$

Under condition (C2), the process $\{N_t\}$ exists for $\tilde \mbp$-a.a. $\tilde\omega.$  For details, see the appendix section below.

Given $\tilde\omega,$ define for $i\in \mathbb Z,$
$$b_i(k)=\left\{\begin{array}{ll}
   \frac{\sum_{j=R-k+1}^R\lambda_i^j}{\mu_i^L}&\text{if } 1\le k\le R,\\
   - \frac{\sum_{j=k-R}^L\mu_i^j}{\mu_i^L}&\text{if } R +1\le k\le R+L-1,
  \end{array}\right.$$
and let \begin{equation*}A_i=\left(
        \begin{array}{cccc}
          0 & 1 & \cdots & 0 \\
           \vdots& \vdots & \ddots &\vdots  \\
          0 & 0 & \cdots & 1 \\
          b_{i}(1) & b_i(2) & \cdots & b_{i}(L+R-1) \\
        \end{array}\right)
\end{equation*} be an $(L+R-1)\times(L+R-1)$ matrix.

Since $A_i$ depends only on $\tilde\omega_i,$ $\{A_i\}_{i\in \mathbb Z}$ is an ergodic sequence of random matrices under $\tilde{\mathbb P}.$ Moreover, under condition (C2), $\tilde{\mathbb E} |\ln \|A_0^{-1}\||+\tilde{\mathbb E} |\ln \|A_0\||<\infty.$ Hence one could use Oseledec's multiplicative ergodic theorem (see \cite{osel}) to the sequence $\{A_i\}_{i\in \mathbb Z}.$ Consequently, we get the Lyapunov
exponents of the sequence $\{A_i\}_{i\in\mathbb Z}$ which we write in increasing order as
$$-\infty<\gamma_1\le \gamma_2\le...\le\gamma_{R+L-1}<\infty.$$
\begin{proposition}[Recurrence/transience criteria]\label{rtcrit} Suppose that condition C holds. Let $\gamma_1\le \gamma_2\le...\le\gamma_{R+L-1}$ be the Lyapunov exponents of the sequence $\{A_i\}_{i\in\mathbb Z}$ under the probability measure $\tilde\mbp.$ Then

\noindent(1) $\gamma_R>0\Rightarrow \tilde P(\lim_{t\rto}N_t=\infty)=1;$

\noindent(2) $\gamma_R=0\Rightarrow  \tilde P(-\infty=\liminf_{t\rto}N_t<\limsup_{t\rto}N_t=\infty)=1;$

\noindent(3) $\gamma_R<0\Rightarrow \tilde P(\lim_{t\rto}N_t=-\infty)=1.$
\end{proposition}
\proof Since the recurrence criteria for $\{N_t\}$ is the same as the embedded process $\{\chi_n\}$ defined below. Proposition \ref{rtcrit} is just a corollary of Theorem A in Letchikov \cite{letc}.\qed

Let $\tau_0=0,$ and for $n\ge 1$ define
$\tau_{n}=\inf\{t>\tau_{n-1}:N_t\neq N_{\tau_{n-1}}\}.$ Since the process $\{N_t\}$ exists, $\tilde P$-a.s., $\tau_n<\infty$ for all $n.$ Note that $\tau_n,\ n\ge 0$ are the consecutive discontinuities of $\{N_t\}.$ Let $\chi_n=N_{\tau_n},\ n\ge 0.$ Then $\{\chi_n\}$ is called the {\it embedded process} of $\{N_t\}.$
Next we study the LLN of  $\{N_t\}.$ Let $T_0=0$ and for $n\ge 1,$ define recursively $$T_{n}=\inf\{t>T_{n-1}:N_t>N_{T_{n-1}}\}.$$  We call $T_n,\ n\ge0$ the ladder times of the process $\{N_t\}.$
Define
\begin{equation}\label{vp}v_{\tilde{\mathbb P}}=\frac{\tilde\mbe \Big( \sum_{r=1}^R\sum_{k\le0}\tilde E_{\theta^{-k}\tilde\omega}\big(\sum_{j=1}^{U_k}\xi_{kj}|N_{T_1}=r\big)\big(\sum_{l=1}^L(-l)\mu_0^l+\sum_{r=1}^Rr\lambda_0^r\big)\Big)}{\sum_{r=1}^R\tilde E(T_1|N_{T_1}=r)}\end{equation}
where $$U_k:=\#\{n:N_{\tau_n}=k,\tau_n<T_1\}$$ is the number of times the embedded process $\{\chi_n\}$ has ever visited $k$ before it hits $[1,\infty),$ and given $\tilde\omega,$ $\xi_{kj},\ k\le 0,\ j\ge 0$ are independent random variables which are also all independent of $U_k$ such that $\tilde P_\omega(\xi_{kj}>t)=e^{-(\sum_{l=1}^L\mu_k^l+\sum_{r=1}^R\lambda_k^r)t},\ t\ge0.$

\begin{theorem}[LLN of $\{N_t\}$]\label{lln}
  Suppose that Condition C holds and $\gamma_R\ge 0.$  Then

\noindent(a) $ \tilde ET_1<\infty\Rightarrow \lim_{t\rto}\frac{N_t}{t}=v_{\tilde\mbp}>0,$ $\tilde P$-a.s.;

\noindent(b) $\tilde ET_1=\infty\Rightarrow \lim_{t\rto}\frac{N_t}{t}=0,$ $\tilde P$-a.s..
\end{theorem}

\begin{remark}\label{rc}
\noindent(i) $\tilde ET_1$ and the velocity $v_{\tilde\mbp}$ of LLN is not given directly in terms of $\tilde\omega.$ They could be calculated  by using the branching structure constructed in \cite{hwb}. Therefore  the formulae of both of $\tilde ET_1$ and $v_{\tilde\mbp}$ could be given explicitly. We treat the case $L=R=2$ to explain the idea in Section \ref{dvp}. See Theorem \ref{vpll} below. The special case $R=1, L> 1$ is discussed in Wang \cite{w}.

\noindent(ii). The case $\gamma_R\le 0$ could be treated in a similar way. We omit this part in this paper.

\noindent(iii). In the proof of part (a) of Theorem \ref{lln}, we only use the  condition
  ``\textbf{(C2')} for some small $\ka>0$ and large $K>0,$ $\tilde\mbp(\lambda_0^1>\ka, \sum_{l=1}^L\mu_0^l+\sum_{r=1}^R\lambda_0^r<K)=1.$"
Under condition (C2'), $\tilde\mbp(\lambda_0^r, \mu_0^l =0 \text{ for certain } 1\le r\le R, 1\le l\le L)>0$ is permitted. So (C2') is weaker than (C2). When proving  part (b) of Theorem \ref{lln} we borrow some results from \cite{brb} where  the uniform ellipticity in (C2) is used.
\end{remark}

To prove Theorem \ref{lln}, we need to prove the LLN of the $h$-skeleton process of $\{N_t\}.$
Fix a number $h>0,$ which will be assumed to be  small enough. For $n\ge 0,$ define $X_n=N_{nh}.$ Then $\{X_n\}$ is a discrete time RWRE with unbounded jumps, which is called the {\it $h$-skeleton process} of $\{N_t\}_{n\ge0}.$ For $i,j\in \mbz$ let $$p_{\tilde\omega}(h,i,j):=\tilde P_{\tilde\omega}(N_h=i+j|N_0=i)$$ be the transition probabilities of $\{X_n\}.$
Define $T_1^h=\inf\{k:X_k> 0\},$ and let
 $$U_k^h=\#\{0\le n<T_1^h:X_n=k\}.$$
\begin{theorem}[LLN of skeleton process]\label{llh}
  Suppose that condition C holds.  Then $\tilde P$-a.s., $\{X_n\}$ is transient to the right, recurrent or transient to the left according as $\gamma_R\ge 0,$ $\gamma_R=0$ or $\gamma_R\le 0.$   Moreover, if $\gamma_R\ge 0,$ then

 \noindent$\tilde ET_1^h=\infty\Rightarrow $ $\tilde P$-a.s., $\lim_{n\rto}\frac{X_n}{n}=0;$

 \noindent$\tilde ET_1^h<\infty\Rightarrow$ $\tilde P$-a.s., $\lim_{n\rto}\frac{X_n}{n}=v_{\tilde\mbp}^h>0,$

\noindent where
$
  v_{\tilde{\mbp}}^h=\frac{\tilde\mbe\z(\sum_{i=1}^\infty\sum_{k\le 0}\tilde E_{\theta^{-k}\tilde\omega}\big(U_k^h|X_{T_1^h=i}\big)\sum_{j\in\mbz} jp_{\tilde\omega}(h, 0,j)\y)}{\sum_{i=1}^\infty \tilde E\big(T_1^h|X_{T_1^h}=i\big)}.
$
\end{theorem}
\begin{remark}
 (i) The $h$-skeleton process $\{X_n\}$ is an RWRE with {\it unbounded jumps}. Under  (C2), we could show that,  for some proper constants $c_0, c_1>0,$ $\tilde\mbp$-a.s., $$p_{\tilde\omega}(h,i,j)<e^{c_0h}e^{-c_1|j|}.$$ Hence, the positive regime of the LLN in Theorem \ref{llh} follows from Theorem \ref{unb}.

 \noindent(ii) One sees easily from Theorem \ref{lln} that $v_{\tilde\mbp}^h$ is indeed independent of $h$ and $v_{\tilde\mbp}^h=v_{\tilde\mbp}.$ For details, see the proof of Theorem \ref{lln} below.
\end{remark}

\section{LLN for RWRE -Proof of Theorem \ref{unb}}\label{pr}

For $n\ge 0$, define $\overline\omega(n)=\theta^{S_n}\omega.$ $\{\overline\omega(n)\}_{n\ge0}$ is an $\Omega^\mathbb N$-valued process. It is usually called ``the environment viewed from particle". Let $$K(\omega,d\omega')=\sum_{j\in\mbz}\omega_{0j}\delta_{\omega'=\theta^j\omega}.$$
\begin{lemma} Under either $P$ or $P_\omega,$ $\{\overline\omega(n)\}_{n\ge0}$ is a Markov chain with transition kernel $K(\omega,\omega').$
\end{lemma}
\proof For test functions $f_1,...,f_n, f_{n+1}$ we have that \begin{eqnarray*}
  &&E_\omega\Big(\prod_{k=1}^{n+1}f_k(\overline\omega(k))\Big)=E_\omega\Big(\prod_{k=1}^{n}f_k(\overline\omega(k))E_\omega^{S_{n}}f_{n+1}(\theta^{S_1}\omega)\Big)\\
&&\quad\quad=E_\omega\Big(\prod_{k=1}^{n}f_k(\overline\omega(k))\sum_{j\in\mbz}\omega_{S_nj}f_{n+1}(\theta^{S_n+j}\omega)\Big)\\
&&\quad\quad=E_\omega\Big(\prod_{k=1}^{n}f_k(\overline\omega(k))\sum_{j\in\mbz}Kf_{n+1}(\overline\omega(n))\Big).
\end{eqnarray*}
Consequently, $\{\overline\omega(n)\}_{n\ge0}$ is a Markov chain under $P_\omega.$ Taking expectation, the above equations also yield the Markov property of  $\{\overline\omega(n)\}_{n\ge0}$ under $P.$ \qed

Whenever $E(T)<\infty,$
define the measures $$ Q(d\omega):=E\z(\displaystyle\sum_{i\ge1}\frac{\mathbf{1}_{S_{T}=i}}{P_\omega(S_{T}=i)}\sum_{k=0}^{T-1}1_{\overline\omega(k)\in d\omega}\y),\ \overline{Q}(d\omega)=\frac{Q(d\omega)}{E(T)}.$$

\begin{lemma}\label{sden}
  Suppose that Condition C holds and $E(T)<\infty.$ Then $Q$ is invariant under the kernel $K,$ that is $$Q(B)=\iint\mathbf1_{\omega'\in B}K(\omega,d\omega')Q(d\omega).$$
Moreover, $Q\sim\mathbb P$ and
$$\frac{dQ}{d\mbp}=\sum_{k\le0}\sum_{i\ge 1}E_{\theta^{-k}\omega}(U_k|S_{T}=i)=:\pi(\omega),$$ where $U_k=\#\{n\le T:S_n=k\}.$
\end{lemma}
\proof By the definition of $Q,$ we have that \begin{eqnarray*}
&&\iint\mathbf1_{\omega'\in B}K(\omega,d\omega')Q(d\omega)=\int P_\omega(\theta^{S_1}\omega\in B)Q(d\omega)\\
&&\quad\quad\quad = E\Big(\sum_{i\ge1}\frac{\mathbf{1}_{S_{T}=i}}{P_\omega(S_{T}=i)}\sum_{k=0}^{T-1}P_{\overline\omega(k)}(\theta^{S_1}\overline\omega(k)\in B)\Big)\\
&&\quad\quad\quad = E\Big(\sum_{i\ge1}\frac{\mathbf{1}_{S_{T}=i}}{P_\omega(S_{T}=i)}\sum_{k=0}^{T-1}P_{\overline\omega(k)}(\overline\omega(k+1)\in B)\Big)\\
&&\quad\quad\quad = E\Big(\sum_{i\ge1}\sum_{k\ge0}P_\omega\z(T\ge k+1,\overline\omega(k+1)\in B\big|S_{T}=i\y) \Big)\\
&&\quad\quad\quad = E\Big(\sum_{i\ge1}\sum_{k\ge1}P_\omega\z(T> k,\overline\omega(k)\in B\big|S_{T}=i\y) \Big)\\
&&\quad\quad\quad \quad\quad\quad+ E\Big(\sum_{i\ge1}P_\omega\z(T<\infty,\theta^i\omega\in B\big|S_{T}=i\y) \Big).\\
\end{eqnarray*}
Since $E(T)<\infty,$ then $P(T<\infty)=P(T> 0)=1.$ This fact together with the stationarity implies that the right-most hand of the last equations equals to
\begin{eqnarray*}
&& E\Big(\sum_{i\ge1}\sum_{k\ge1}P_\omega\z(T> k,\overline\omega(k)\in B\big|S_{T}=i\y) \Big)\\
&&\quad\quad\quad\quad\quad+ E\Big(\sum_{i\ge1}P_\omega\z(T> 0,\overline\omega(0)\in B\big|S_{T}=i\y) \Big)\\
&&\quad\quad=E\Big(\sum_{i\ge0}\sum_{k\ge0}P_\omega\z(T> k,\overline\omega(k)\in B\big|S_{T}=i\y) \Big)=Q(B).
\end{eqnarray*}
The first part of the lemma follows.
To prove the second part, for testing function $f(\omega),$ we have that
\begin{eqnarray}
 && \int f dQ=E\Big(\sum_{i\ge1}\frac{\mathbf{1}_{S_{T}=i}}{P_\omega(S_{T}=i)}\sum_{k=0}^{T-1}f(\overline\omega(k))\Big)\no\\
&&\quad\quad\quad\ =E\Big(\sum_{i\ge1}\frac{\mathbf{1}_{S_{T}=i}}{P_\omega(S_{T}=i)}\sum_{k=0}^{T-1}f(\theta^{S_k}\omega)\Big)\nonumber\\
&&\quad\quad\quad\  =E\Big(\sum_{i\ge1}\frac{\mathbf{1}_{S_{T}=i}}{P_\omega(S_{T}=i)}\sum_{k\le 0}U_kf(\theta^{k}\omega)\Big)\no\\
&&\quad\quad\quad\ =\mbe\Big(\sum_{i\ge1}E_\omega\Big(\sum_{k\le 0}U_kf(\theta^{k}\omega)\big|S_{T}=i\Big)\Big)\nonumber\\
&&\quad\quad\quad\ =\mbe\Big(f(\omega)\sum_{i\ge1}\sum_{k\le 0}E_{\theta^{-k}\omega}\Big(U_k\big|S_{T}=i\Big)\Big).\label{eqf}
\end{eqnarray}
By stationarity, we have that
\begin{eqnarray*}
\mbe\Big(\sum_{i\ge1}\sum_{k\le 0}E_{\theta^{-k}\omega}\Big(U_k\big|S_{T}=i\Big)\Big)=\sum_{i\ge1}\sum_{k\le 0}E\Big(U_k\big|S_{T}=i\Big)\le E(T)<\infty.
\end{eqnarray*}
Then it follows from  (\ref{eqf})  that $Q\sim \mbp$ and $$\frac{dQ}{d\mbp}=\sum_{k\le0}\sum_{i\ge 1}E_{\theta^{-k}\omega}(U_k|S_{T}=i).$$
\qed

\begin{lemma}\label{sbo}
  Under the conditions of Lemma \ref{sden}, $\{\overline\omega(n)\}$ is stationary and ergodic under the probability measure $\overline Q\times P_\omega.$
\end{lemma}
With Lemma \ref{sden} in hand, Lemma \ref{sbo} follows similarly as  Sznitman \cite{S02}, Theorem 1.2 or  Zeitouni \cite{ze04}, Corollary 2.1.25.

Define the local drift $d(x,\omega)=E_\omega^x(S_1-S_0)$ and set $$M_n=S_n-S_0-\sum_{k=0}^{n-1}d(S_n,\omega).$$
\begin{lemma}\label{coy}
  Under $P_\omega,$ $\{M_n\}$ is a martingale and $P$-a.s., $\lim_{n\rto}\frac{M_n}{n}=0.$
\end{lemma}
\proof Note that \begin{equation*}\begin{split}E_\omega&(M_n-M_{n-1}|M_{n-1},...,M_{0})\\
&=E_\omega(S_n-S_{n-1}-d(S_{n-1},\omega)|S_{n-1},...,S_0)=0.\end{split}\end{equation*} Then $\{M_n\}$ is a martingale under $P_\omega.$ By (B3), there exist some constants $c_2>0$ and $0<\varepsilon_1<\varepsilon_0$ such that $E_\omega(|M_n-M_{n-1}|^{2+\varepsilon_1})<c_2.$   Then one follows from Theorem 3.2 in \cite{lv} that there exists constant $c_3>0$ such that for $\lambda>0,$ and $n$ large enough,
\begin{equation}\label{azm}P(|M_n|>\sqrt{n}\lambda)\le \frac{c_3}{\lambda^{2+\varepsilon_1}}.\end{equation}
Choosing  $0<\varepsilon_2<\frac{1}{2}$ properly, we have $(2+\varepsilon_{1})\varepsilon_2>1.$
Setting $\lambda=n^{\varepsilon_2}$ in (\ref{azm}), it follows that
  $$P(|M_n|>n^{\frac{1}{2}+\varepsilon_2})\le c_3n^{-(2+\varepsilon_1)\varepsilon_2}.$$
An application of Borel-Cantelli's lemma yields that $$P\text{-a.s., } \lim_{n\rto}\frac{M_n}{n}=0.$$ \qed

\noindent{\it Proof of Theorem \ref{unb}:}
By Lemma \ref{sbo}, $\{\overline\omega(n)\}$ is a stationary and ergodic sequence under the measure $\overline Q\times P_\omega.$ Using Birkhoff's ergodic theorem, we have that for $\overline Q$-a.a. or $\mbp$-a.a. $\omega,$ $P_\omega$-a.s.,
\begin{equation}\label{cdo}
  \lim_{n\rto}\frac{1}{n}\sum_{k=0}^{n-1}d(S_k,\omega)=\lim_{n\rto}\frac{1}{n}\sum_{k=0}^{n-1}d(0,\overline\omega(k))=\int d(0,\omega)d\overline Q.
\end{equation}
We conclude from Lemma \ref{coy} and (\ref{cdo}) that $P$-a.s.,
$$\lim_{n\rto}\frac{S_n}{n}=\int d(0,\omega)d\overline Q=:v_\mbp.$$
It follows from Lemma \ref{sden} that
\begin{equation*}
  \begin{split}
  v_\mbp=\int d(0,\omega)d\overline Q=\frac{\mbe\z(\sum_{i=1}^\infty\sum_{k\le 0}E_{\theta^{-k}\omega}(U_k|S_{T=i})\sum_{j\in\mbz} j\omega_{0j}\y)}{\sum_{i=1}^\infty E(T|S_{T}=i)}.
\end{split}
\end{equation*}
Theorem \ref{unb} is proved. \qed

\section{LLN for BDPRE with bounded jumps}\label{pec}
\subsection{LLN for the $h$-skeleton process-Proof of Theorem \ref{llh}}\label{sllh}
The recurrence criteria for $\{X_n\}$ in Theorem \ref{llh} follows directly from Proposition \ref{rtcrit}. The zero speed regime of the LLN for $\{X_n\}$ follows from the counterpart of Theorem \ref{lln}. Next we prove
the non-zero speed regime of the LLN in Theorem \ref{llh}. To begin with, we estimate the tail probability of the transition probability
$p_{\tilde\omega}(h,i,j):=\tilde P_{\tilde\omega}(N_h=i+j|N_0=i)$
of $\{X_n\}.$
The following lemma shows that for fixed $i\in \mathbb Z,$ $\tilde{\mathbb P}$-a.s.,  $p_{\tilde\omega}(h,i,j)$ decays exponentially to $0$ as $|j|\rto.$
\begin{lemma}\label{esp}
  Suppose that Condition (C2) is satisfied. Then for $\tilde\mbp$-a.a. $\tilde\omega,$ there exist $0<c_0<\infty$ and $0<c_1<\infty,$ which are independent of $h$ and $\tilde\omega,$ such that for $|j|> \max\{L,R\},$ \begin{equation}\label{epij}
  p_{\tilde\omega}(h,i,j)<e^{c_0h}e^{-c_1|j|},
\end{equation} where $c_1$ could be made arbitrarily large by adjusting the value of $c_0.$
\end{lemma}
\proof We prove only the case $j>\max\{L,R\}.$ The case $j<-\max\{L,R\}$ follows similarly. Let $m=[\frac{j}{R}].$ Since at each discontinuity, the process $\{N_t\}$ jumps at most a distance $R$ to the right, then starting from $i,$ in order to reach $i+j,$  $\{N_t\}$ has at least $m$ discontinuities in a time period of length $h.$ Let $\eta_k,\ k=1,...,m$ be these discontinuities and $\tau_k$ be the waiting time after $\eta_k$ until the process $\{N_t\}$ leaves $N_{\eta_k}.$ Then we have that
\begin{equation}\label{ptau}
  p_{\tilde\omega}(h,i,j)= \tilde P_{\tilde\omega}(N_h=i+j|N_0=i)\le \tilde P^i_{\tilde\omega}(\tau_1+...+\tau_m\le h).
\end{equation}
Note that  $\tilde P_{\tilde\omega}^i(\tau_k>t|N_{\eta_k}=i_k)=e^{-s_kt}$ for $t>0$ where $s_k=\big(\sum_{l=1}^L\mu_{i_k}^l+\sum_{r=1}^R\lambda_{i_k}^r\big).$ Moreover under $\tilde P_{\tilde\omega}^i,$ $\tau_k,\ k=1,...,m$ are mutually independent. Then by Chebycheff's bound, it follows that  for $\lambda<0,$ $K=(L+R)M,$ $\ka=(L+R)\varepsilon,$
 \begin{eqnarray*}
&&\tilde P^i_{\tilde\omega}(\tau_1+...+\tau_m\le h)\\
&&\quad=\sum_{i_1,...,i_k}\tilde P^i_{\tilde\omega}\Big(\sum_{k=1}^m\tau_k\le h\Big|N_{\eta_k}=i_k,1\le k\le m\Big)\tilde P^i_{\tilde\omega}(N_{\eta_k}=i_k,1\le k\le m)\\
&&\quad\le\sum_{i_1,...,i_k}e^{-\lambda h}\tilde E_{\tilde\omega}^i\Big(e^{\lambda\sum_{k=1}^m\tau_k}\Big|N_{\eta_k}=i_k,1\le k\le m\Big)\\
&&\quad\quad\quad\quad\quad\quad\quad\quad\quad\quad\quad\quad\quad\times\tilde P_{\tilde\omega}^i(N_{\eta_k}=i_k,1\le k\le m)\\
&&\quad= e^{-\lambda h}\sum_{i_1,...,i_k}\prod_{k=1}^m \frac{s_k}{s_k-\lambda}\tilde P^i_{\tilde\omega}(N_{\eta_k}=i_k,1\le k\le m)\\
&&\quad\le  e^{-\lambda h}\Big(\frac{K}{\ka-\lambda}\Big)^m, \ \tilde\mbp\text{-a.s.,}
\end{eqnarray*}
where the last inequality follows from condition (C2). Substituting the above estimation to (\ref{ptau}), we have that $\mbp$-a.s.,
$$p_{\tilde\omega}(h,i,j)\le e^{-\lambda h}\Big(\frac{K}{\ka-\lambda}\Big)^m\le e^{-\lambda h}e^{\frac{j}{R}(\log K-\log(\ka-\lambda))}.$$
By choosing $\overline{\lambda}<0$ properly and letting  $c_0=-\overline{\lambda},$  $c_1=(\log(\ka-\overline\lambda)-\log K)/R,$ (\ref{epij}) is proved. Of course, we could make $c_1$ arbitrarily large by adjusting the value of $\overline\lambda.$ \qed

\begin{lemma}\label{sep} Under the measure $\tilde\mbp,$
  $\{p_{\tilde\omega}(h,i,j)\}_{i\in\mbz}$ is a stationary and ergodic sequence.
\end{lemma}
For the proof of the lemma, refer to Durrett \cite{dur}.
\qed

One follows from condition (C2) that for some small constant $c_4>0$ \begin{equation}\label{ell}p_{\tilde\omega}(h,0,1)\ge c_4.\end{equation}
Taking Lemma \ref{esp}, Lemma \ref{sep} and (\ref{ell}) together,  we could use Theorem \ref{unb} to conclude that
  $$\tilde E(T_1^h)<\infty\Rightarrow \tilde{P}\text{-a.s., } \lim_{n\rto}\frac{X_n}{n}=v_{\tilde\mbp}^h.$$
  We complete the proof of Theorem \ref{llh}. \qed

\subsection{The  LLN of $\{N_t\}$-Proof of Theorem \ref{lln}}
In this subsection,  using the LLN (nonzero speed regime) of the $h$-skeleton process proved in Subsection  \ref{sllh}, we prove the LLN of the process $\{N_t\}.$ Firstly, we show that for $h$ small enough, $\tilde E(T_1^h)<\infty$ whenever $\tilde E(T_1)<\infty.$ Therefore, under the condition of part (a) of Theorem \ref{lln}, the condition of the nonzero speed regime of LLN in Theorem \ref{llh} is also satisfied for $h$ small enough.
\begin{lemma}\label{tht}
 Fix $s>1.$ Suppose that condition (C2) holds.   Then, for $n$ large enough, $\tilde P$-a.s., $T_1^{\frac{1}{n^{s}}}\frac{1}{n^s}\le T_1+\frac{1}{n^s}.$
Consequently, we have that for $n$ large, if $\tilde E(T_1)<\infty,$ then $\tilde E\Big(T_1^{\frac{1}{n^{s}}}\Big)<\infty$ and $\tilde P$-a.s., $\lim_{n\rto}T_1^{\frac{1}{n^{s}}}\frac{1}{n^s}=T_1.$
\end{lemma}
\proof Given ${\tilde\omega},$ we have that
\begin{eqnarray}\label{tns}
  &&\tilde P_{\tilde\omega}\Big(T_1^{\frac{1}{n^{s}}}\frac{1}{n^s}> T_1+\frac{1}{n^s}\Big)\no\\
  &&\quad\quad\quad\le \tilde P_{\tilde\omega}\Big(\{N_t\} \text{ has at least one jump in } [T_1,T_1+\frac{1}{n^s})\Big)\nonumber\\
&&\quad\quad\quad \le \tilde P_{\tilde\omega}\Big(\{N_t\} \text{ leaves  }N_{T_1} \text{ within time } \frac{1}{n^s}\Big)\nonumber\\
&&\quad\quad\quad =\sum_{r=1}^R\tilde P_{\tilde\omega}(N_{T_1}=r)\tilde P_{\tilde\omega}\Big(\nu_r<\frac{1}{n^s}\big|N_{T_1}=r\Big)
\end{eqnarray}
where $\nu_r$ is the waiting time at state $r$ until the next jump of $\{N_t\}$ happens.
Since $\nu_r$ is exponentially distributed with parameter $\zeta_r:=\sum_{l=1}^L\mu_r^l+\sum_{k=1}^R\lambda_r^k.$
Using condition (C2), with $K=(L+R)M,$ it follows that $\tilde\mbp$-a.s., the right-most hand of  (\ref{tns}) equals to
$$\sum_{r=1}^R\tilde P_{\tilde\omega}(N_{T_1}=r)(1-e^{-\zeta_r \frac{1}{n^s}})\le 1-e^{-K \frac{1}{n^s}}. $$
Noting that $s>1,$ then we have that
$$\sum_{n=1}^\infty \tilde P\Big(T_1^{\frac{1}{n^{s}}}\frac{1}{n^s}> T_1+\frac{1}{n^s}\Big)<\infty.$$
An application of Borel Cantelli's lemma yields that for $n$ large, $\tilde P$-a.s.,
\begin{equation}\label{tst}
  T_1^{\frac{1}{n^{s}}}\frac{1}{n^s}\le T_1+\frac{1}{n^s}.
\end{equation}
Since $\tilde P$-a.s., $T_1^hh>T_1,$ then one follows from (\ref{tst}) that $\tilde P$-a.s.,
$\lim_{n\rto}T_1^{\frac{1}{n^{s}}}\frac{1}{n^s}= T_1.$ \qed

\noindent{\it Proof of Theorem \ref{lln}:} We prove part (a) first. Suppose $\gamma_R\ge0$ and $\tilde E(T_1)<\infty.$ Fix $h>0$ small enough. For any $t>0,$ there is a unique number $n_t$ such that $n_th\le t<(n_{t}+1)h.$ Let $J_t$ be the number of jumps of $\{N_t\}$ in the time interval $[n_th, (n_{t}+1)h).$ Then, under condition (C2), a similar argument as the proof of Lemma \ref{esp} yields that there exist some positive constants $c_5$ and $c_6$ such that
\begin{equation}\label{bpk}\tilde P(J_t>n)\le e^{c_5h}e^{-nc_6}.\end{equation}
Note that the bound in (\ref{bpk}) is independent of $t.$ Applying Borel-Cantelli lemma, we have that $\tilde P$-a.s., \begin{equation}\label{jt0}\lim_{n\rto}\frac{J_t}{n}=0 \end{equation} uniformly in $t.$
Since at each discontinuity, $\{N_t\}$ jumps at most a distance $L$ to the left or at most a distance $R$ to the right, we have that
\begin{equation}\label{jbt}
  \frac{N_{n_th}-J_tL}{(n_t+1)h}\le \frac{N_t}{t}\le  \frac{N_{n_th}+J_tR}{n_th}.
\end{equation}
For $h$ small enough, by Lemma \ref{tht}, $\tilde E(T_1^h)<\infty$. Then we have from the nonzero speed regime of LLN in Theorem \ref{llh}, (\ref{jt0}) and (\ref{jbt}) that $\tilde P$-a.s.,
\begin{equation}\label{lmh}
  \lim_{t\rto}\frac{N_t}{t}=\frac{v_{\tilde\mbp}^h}{h}.
\end{equation}
To finish the proof of part (a) in Theorem \ref{lln}, it suffices to show that, for all $h>0,$
\begin{equation}\label{vph}\frac{v_{\tilde\mbp}^h}{h}=v_{\tilde\mbp}.\end{equation}
Indeed, since  $\frac{N_t}{t}$ is independent of $h$ and the limit in (\ref{lmh}) exists, then $\frac{v_{\tilde\mbp}^h}{h}$ is independent of $h.$
Therefore we have that
\begin{eqnarray}\label{lhd}
&&  \frac{v_{\tilde\mbp}^h}{h}=\lim_{h\rightarrow 0}\frac{v^h_{\tilde\mbp}}{h}=\lim_{h\rightarrow0}\frac{\tilde\mbe\z(\sum_{i=1}^\infty\sum_{k\le 0}\tilde E_{\theta^{-k}{\tilde\omega}}(U_k^h|X_{T_1^h=i})\sum_{j\in\mbz} jp_{\tilde\omega}(h, 0,j)\y)}{h\sum_{i=1}^\infty \tilde E(T_1^h|X_{T_1^h}=i)}\nonumber\\
&&\quad\quad\quad=\lim_{h\rightarrow0}\frac{\tilde\mbe\z(\sum_{i=1}^\infty\sum_{k\le 0}\tilde E_{\theta^{-k}{\tilde\omega}}(hU_k^h|N_{hT_1^h=i})\sum_{j\in\mbz} j\frac{p_{\tilde\omega}(h, 0,j)}{h}\y)}{\sum_{i=1}^\infty \tilde E(h T_1^h|N_{hT_1^h}=i)}.
\end{eqnarray}
Note that by Lemma \ref{esp}, under condition (C2), $\tilde P$-a.s., $\sum_{j\in\mbz} j\frac{p_{\tilde\omega}(h, 0,j)}{h}$ is uniformly bounded from above.
Moreover, by Lemma \ref{tht}, $\tilde P$-a.s., $hT^h_1$ is bounded by $T_1+h$ for $h$ small, and by stationarity
\begin{equation*}\begin{split}\tilde\mbe\Big(\sum_{i=1}^\infty\sum_{k\le 0}&\tilde E_{\theta^{-k}{\tilde\omega}}(hU_k^h|N_{hT_1^h=i})\Big)=\sum_{i=1}^\infty\sum_{k\le 0}\tilde\mbe\z(\tilde E_{\theta^{-k}{\tilde\omega}}(hU_k^h|N_{hT_1^h=i})\y)\\
&=\sum_{i=1}^\infty\sum_{k\le 0}\tilde E\z(hU_k^h|N_{hT_1^h=i}\y)=\sum_{i=1}^\infty \tilde E\z(hT_1^h|N_{hT_1^h=i}\y).\end{split}\end{equation*}
Since $\tilde P$-a.s., $\lim_{h\rto 0}hT^h_1=T_1,$ then by the above discussion, the condition of dominated convergence theorem is satisfied.
Note also that $\tilde P$-a.s.,
\begin{equation}\label{uhk}\lim_{h\rightarrow 0}hU^h_k=\lim_{h\rightarrow0}\sum_{j=0}^{T_1^h-1}h\mathbf1_{N_{jh}=k}=\int_{0}^{T_1}\mathbf1_{N_t=k}dt\overset{\mathcal D}{=}\sum_{i=1}^{U_k}\xi_{ki},\end{equation}
where $ ``A\overset{\mathcal D}{=}B"$ means that A equals to B in $\tilde P_{\tilde\omega}$ distribution.
Since $\tilde P$-a.s., $\lim_{h\rightarrow 0}hT^h_1=T_1,$  $hT^h_1\ge T_1,$ and the paths of $\{N_t\}$ are right continuous,  then by dominated convergence, it follows from (\ref{lhd}) and (\ref{uhk}) that
\begin{eqnarray*}
\frac{v_{\tilde\mbp}^h}{h}&=&\lim_{h\rightarrow0}\frac{\tilde\mbe\z(\sum_{i=1}^\infty\sum_{k\le 0}\tilde E_{\theta^{-k}{\tilde\omega}}(hU_k^h|N_{hT_1^h=i})\sum_{j\in\mbz} j\frac{p_{\tilde\omega}(h, 0,j)}{h}\y)}{\sum_{i=1}^\infty\tilde E(h T_1^h|N_{hT_1^h}=i)}\\
&=&\frac{\tilde\mbe \Big( \sum_{r=1}^R\sum_{k\le0}\tilde E_{\theta^{-k}{\tilde\omega}}\big(\sum_{j=1}^{U_k}\xi_{kj}|N_{T_1}=r\big)\big(\sum_{l=1}^L(-l)\mu_0^l+\sum_{r=1}^Rr\lambda_0^r\big)\Big)}{\sum_{r=1}^R\tilde E(T_1|N_{T_1}=r)}\\
&=&v_{\tilde\mbp}>0,
\end{eqnarray*}
where in the second equality, we use the facts $\tilde P(N_{T_1}=r)=0$ for $r>R$ and
$$\lim_{h\rightarrow 0}\frac{p_{\tilde\omega}(h,0,j)}{h}=\z\{\begin{array}{ll}
                                         \lambda_0^j & \text{if } j=1,...,R, \\
                                         \mu_0^j & \text{if } j=1,...,L, \\
                                         0 & \text{if }j<-L \text{ and } j>R.
                                       \end{array}\y.
$$
 Consequently, (\ref{vph}) is proved and part (a) of Theorem \ref{lln} follows.

Next, we turn to prove part (b) of Theorem \ref{lln}. Suppose $\gamma_R\ge 0$ and $\tilde E(T_1)=\infty.$ Recall that $\tau_n,\ n\ge 0$ are the consecutive discontinuities of $\{N_t\}$ and $\{\chi_n\}=\{N_{\tau_n}\}$ is the embedded process.  Let $\overline T_1:=\inf\{k:\chi_k>0\}.$
As the first step, we show that
\begin{equation}\label{tt}
 \tilde E(\overline T_1)=\infty.
\end{equation}
For this purpose, for $k\le 0$ let $\overline{U}_k=\#\{n<\overline T_1:\chi_n=k\}.$ Then Wald's equation implies that,  $$\tilde E_{\tilde\omega}(T_1)=\sum_{k\le 0}\tilde E_{\tilde\omega}\Big(\sum_{i=1}^{\overline U_k}\xi_{ki}\Big)=\sum_{k\le 0}\tilde E_{\tilde\omega}(\overline U_k)\tilde E_{\tilde\omega}(\xi_{k1}).$$
But condition (C2) implies that $\tilde\mbp$-a.s., $$\frac{1}{K}<\tilde E_{\tilde\omega}(\xi_{k1})=\frac{1}{\sum_{l=1}^L\mu_k^l+\sum_{r=1}^R\lambda_k^r}<\frac{1}{\ka}$$
with $K=(L+R)M$ and $\ka=(L+R)\varepsilon.$
Then it follows  that $\tilde\mbp$-a.s.,
$$\frac{1}{K}\sum_{k\le 0}\tilde E_{\tilde\omega}(\overline U_k)\le \tilde E_{\tilde\omega}(T_1)\le \frac{1}{\ka}\sum_{k\le 0}\tilde E_{\tilde\omega}(\overline U_k).$$
Taking expectation, we have that
$$\frac{\tilde E(\overline T_1)}{K}\le\tilde E(T_1)\le \frac{\tilde E(\overline T_1)}{\ka}.$$
Therefore (\ref{tt}) follows. Then it follows from Br\'emont \cite{brb} (See, Proposition 9.1, Theorem 9.2 and Corollary 9.3 therein.) that $\tilde P$-a.s., \begin{equation}\label{limnt}\lim_{n\rto}\frac{\chi_n}{n}=\lim_{n\rto}\frac{N_{\tau_n}}{n}=0.\end{equation}
As the second step, we show that for $h>0$ small enough and $n$ large enough
\begin{equation}\label{tnh}\tilde P\text{-a.s., } \frac{\tau_n}{n}>h.\end{equation}
For $k\ge 1,$ let $\nu_k=\tau_k-\tau_{k-1}.$ Then under $\tilde P_{\tilde\omega},$ $\nu_k,\ k\ge 1$ are mutually independent and all exponentially distributed.
By Chebycheff's bound, it follows that for $\lambda<0,$ \begin{eqnarray*}\label{lde}
  &&\tilde P_{\tilde\omega}(\tau_n\le nh)\le e^{-\lambda nh}\tilde E_{\tilde\omega}(e^{\lambda \sum_{k=1}^n\nu_k})=e^{-\lambda nh}\prod_{k=1}^n \tilde E_{\tilde\omega}(e^{\lambda \nu_k})\nonumber\\
&&\quad\quad\quad\quad\le e^{-n(\lambda h-\log K+\log(\ka-\lambda))}
\end{eqnarray*}
with $\ka=(L+R)\varepsilon$ and $K=(L+R)M.$
Choosing properly $\overline\lambda<0$ and $h$ small enough, we could make $c(h):=\overline\lambda h-\log K+\log(\ka-\overline\lambda)$ a strictly positive number.
Consequently we have $\sum_{n=1}^\infty \tilde P_{\tilde\omega}(\tau_n\le nh)<\infty.$ By applying Borel-Cantelli lemma, we have (\ref{tnh}).

Finally, (\ref{limnt}) and (\ref{tnh}) imply that $\tilde P$-a.s., \begin{equation}\label{lmntt}\lim_{n\rto}\frac{N_{\tau_n}}{\tau_n}=0.\end{equation}
For $t>0,$ there is a unique random number $n_t$ such that $\tau_{n_t}\le t<\tau_{n_t+1}.$
Since at each discontinuity, $\{N_t\}$ jumps at most a distance $L$ to the left or at most a distance $R$ to the right, we have  \begin{equation}\label{jbnt}\frac{N_{\tau_{n_t}}-L}{\tau_{n_t+1}}\le \frac{N_t}{t}\le \frac{N_{\tau_{n_t}}+R}{\tau_{n_t}}\end{equation}
It follows from (\ref{lmntt}) and (\ref{jbnt}) that $\tilde P$-a.s., $$\lim_{t\rto}\frac{N_t}{t}=0.$$
Part (b) of Theorem \ref{lln} is proved.
\qed

\section{Discussion of the velocity $v_{\tilde\mbp}$}\label{dvp}
In this section, we let $R=L=2$ and discuss the asymptotic velocity $v_{\tilde\mbp}.$ We have
\begin{theorem}\label{vpll}
  Let $\pi({\tilde\omega})$ and $D({\tilde\omega})$ be as in (\ref{do}) and (\ref{pio}) below. Suppose $L=R=2$ and $\tilde\mbe(\pi({\tilde\omega}))<\infty.$ Then
$\tilde P$-a.s., $$\lim_{t\rto}\frac{N_t}{t}=\frac{\tilde\mbe\z(\pi({\tilde\omega})(2\lambda_0^2+\lambda_0^1-\mu_0^1-2\mu_0^2)\y)}{\tilde\mbe(D({\tilde\omega}))}.$$
\end{theorem}
\proof Theorem \ref{vpll} is just a special case of Theorem \ref{lln}. We need only to calculate $v_{\tilde\mbp}.$
Recall that $\overline T_1=\inf\{n\ge 0:\chi_n > 0\}$ is the first ladder time
of the embedded process. Define $$\overline U_k=\#\{n\le \overline T_1:\chi_n=k\}$$ which is the occupation time at  state $k$ of the embedded process before $\overline{T}_1.$  In \cite{hwb}, the authors show that $\overline{U}_k$ could be written as the functional of a multitype branching process $\{Z_n\}_{n\le 1}$ whose offspring matrices are the function of the environment ${\tilde\omega}.$  To introduce those results, we  need to introduce some notations.

Fix $a<b.$ Let $\partial^+[a,b]=\{b,b+1\}$ and
$\partial^-[a,b]=\{a,a-1\}$ be the positive and negative boundaries
of $[a,b]$ correspondingly. For $k\in (a,b),$
$\zeta\in\partial^+[a,b]\cup\partial^-[a,b],$ define
\begin{equation}\label{epk}\no
\mcp_k(a,b,\zeta)=\tilde P_{\tilde\omega}^k(\{\chi_n\}\text{ exits the interval
}[a+1,b-1]\text{ at }\zeta).\end{equation}
For $j=1,2$ let $p_i^j=\frac{\lambda_i^j}{\mu_i^1+\mu_i^2+\lambda_i^1+\lambda_i^2}$ and $q_i^j=\frac{\mu_i^j}{\mu_i^1+\mu_i^2+\lambda_i^1+\lambda_i^2}.$

 Writing
$\mcp_k(a,b,\zeta)$ as $\mcp_k(\zeta)$ temporarily, one follows from
Markov property that
\begin{equation*}
  \mcp_k(\zeta)=p_k^2\mcp_{k+2}(\zeta)+p_k^1\mcp_{k+1}(\zeta)+q_k^1\mcp_{k-1}(\zeta)+q_k^2\mcp_{k-2}(\zeta),
\end{equation*}
which leads to the following matrix form
$$  V_k(\zeta)=M_kV_{k+1}(\zeta)$$  where
$$  V_k(\zeta):=\left(\begin{array}{c}
             (\mcp_{k-1}-\mcp_{k-2})(\zeta) \\
             (\mcp_{k}-\mcp_{k-1})(\zeta) \\
             (\mcp_{k+1}-\mcp_{k})(\zeta) \\
           \end{array}
         \right), M_k=\left(
    \begin{array}{ccc}
      -\frac{\mu_k^1+\mu_k^2}{\mu_k^2} & \frac{\lambda_k^1+\lambda_k^2}{\mu_k^2} & \frac{\lambda_k^2}{\mu_k^2} \\
     1  &  0& 0\\
      0& 1 & 0 \\
    \end{array}
  \right).
$$
Then  $\mcp_k(a,b,b)$ and  $\mcp_k(a,b,b+1)$ follow from some standard procedure.
They could be expressed in terms of $\{M_i\}_{i\in\mathbb{Z}}.$
For $k\le i,$ let
$$f_k(i,i+1)=\mcp_k(-\infty,i+1,i+1)\text{ and } f_k(i,i+2)=\mcp_k(-\infty,i+1,i+2).$$
For $i\in\mbz,$ let
\begin{eqnarray*}
&&\alpha_{i,1}=\frac{q_i^1p_{i-1}^1}{1-q_{i-1}^1f_{i-2}(i-2,i-1)-q_{i-1}^2f_{i-3}(i-2,i-1)},\\
&&\alpha_{i,3}=\frac{q_i^1p_{i-1}^2}{1-q_{i-1}^1f_{i-2}(i-2,i-1)-q_{i-1}^2f_{i-3}(i-2,i-1)},\\
&&\beta_{i,1}=\frac{q_i^2f_{i-2}(i-2,i-1)p_{i-1}^1}{1-q_{i-1}^1f_{i-2}(i-2,i-1)-q_{i-1}^2f_{i-3}(i-2,i-1)},\\
&&\beta_{i,3}=\frac{q_i^2f_{i-2}(i-2,i-1)p_{i-1}^2}{1-q_{i-1}^1f_{i-2}(i-2,i-1)-q_{i-1}^2f_{i-3}(i-2,i-1)},\\
&&\gamma_{i,1}=\frac{q_{i+1}^2p_{i-1}^1}{1-q_{i-1}^1f_{i-2}(i-2,i-1)-q_{i-1}^2f_{i-3}(i-2,i-1)},\\
&&\gamma_{i,3}=\frac{q_{i+1}^2p_{i-1}^2}{1-q_{i-1}^1f_{i-2}(i-2,i-1)-q_{i-1}^2f_{i-3}(i-2,i-1)},\\
&&\ \alpha_{i,2}:=q_i^1-\alpha_{i,1}-\alpha_{i,3};\ \beta_{i,2}:=q_i^2-\beta_{i,1}-\beta_{i,3};\ \gamma_{i,2}:=q_{i+1}^2-\gamma_{i,1}-\gamma_{i,3}.
  \end{eqnarray*}
Set \begin{equation}\label{du}
u_1:=\z(\frac{\alpha_{1,1}}{\alpha_{1,1}+\alpha_{1,2}+\alpha_{1,3}},\frac{\alpha_{1,2}}{\alpha_{1,1}+\alpha_{1,2}+\alpha_{1,3}},\frac{\alpha_{1,3}}{\alpha_{1,1}+\alpha_{1,2}+\alpha_{1,3}},0,...,0\y)\in \mathbb R^9,
\end{equation}
and for $i\le 0$ define
$x_i=\frac{\alpha_{i,1}}{1-\alpha_{i,1}-\alpha_{i,2}-\beta_{i,1}-\beta_{i,2}},$
$y_i=\frac{\alpha_{i,2}}{1-\alpha_{i,1}-\alpha_{i,2}-\beta_{i,1}-\beta_{i,2}},$
$z_i=\frac{\beta_{i,1}}{1-\alpha_{i,1}-\alpha_{i,2}-\beta_{i,1}-\beta_{i,2}},$
$w_i=\frac{\beta_{i,2}}{1-\alpha_{i,1}-\alpha_{i,2}-\beta_{i,1}-\beta_{i,2}},$
 $1-v_i=\frac{\gamma_{i,3}}{\beta_{i+1,2}},$
 $s_i=\frac{\alpha_{i,3}}{\alpha_{i,3}+\beta_{i,3}}$ and
 $t_i=\frac{\gamma_{i,1}}{\gamma_{i,1}+\gamma_{i,2}}.$ Define also the matrices
\begin{equation}\label{q}
  Q_i=\left(
    \begin{array}{ccccccccc}
      x_i & y_i & 0 & z_i & w_i & 0 & 0 & 0 & 0 \\
       x_i & y_i & s_i & z_i & w_i & 1-s_i & 0 & 0 & 0 \\
       x_i & y_i & 0 & z_i & w_i & 0 & 0 & 0 & 0 \\
      x_i & y_i & 0 & z_i & w_i & 0 & t_i & 1-t_i  & 0 \\
      x_iv_i & y_iv_i & s_iv_i & z_iv_i & w_iv_i & (1-s_i)v_i & t_iv_i & (1-t_i) v_i & 1-v_i \\
    x_i & y_i & 0 & z_i & w_i & 0 & t_i & 1-t_i  & 0\\
      x_i & y_i & 0 & z_i & w_i & 0 & 0 & 0 & 0 \\
       x_i & y_i & s_i & z_i & w_i & 1-s_i & 0 & 0 & 0 \\
      x_i & y_i & 0 & z_i & w_i & 0 & 0 & 0 & 0 \\
    \end{array}
  \right).
\end{equation}
We have the following proposition which could be found in \cite{hwb}.
\begin{proposition} Suppose that $\gamma_R\ge 0.$ For $\tilde\mbp$-a.a. ${\tilde\omega},$ there exists a 9-type branching process $\{Z_n\},$ whose mean offspring matrices are as in (\ref{q}) and initial distribution mean is $u_1$ defined in (\ref{du}), such that in $\tilde P_{\tilde\omega}$-distribution, with $\mathbf v_1=(1,1,1,0,0,0,1,1,1)^T$ and $\mathbf v_2=(1,1,0,1,1,0,1,1,0)^T,$
  $$\overline U_k=Z_{k+1}\mathbf v_1+Z_k\mathbf v_2, \ k\le0.$$
Moreover, with the empty product being identity,
 \begin{equation}\label{enif}\begin{split}
&\tilde E_{\tilde\omega}(\overline U_k|\chi_{\overline T_1}=2)+\tilde E_{\tilde\omega}(\overline U_k|\chi_{\overline T_1}=1)\\
 &=\Big(\frac{\alpha_{1,1}}{\alpha_{1,1}+\alpha_{1,2}},\frac{\alpha_{1,2}}{\alpha_{1,1}+\alpha_{1,2}},1,0,...,0\Big)(Q_0Q_{-1}\cdots
Q_{k+1}\mathbf v_1+Q_0Q_{-1}\cdots Q_{k}\mathbf v_2).
\end{split}
  \end{equation}
\end{proposition}
Since after the $i$th visit of state $k,$  $\{N_t\}$  will wait here an exponentially distributed  time $\xi_{ki}$ with parameter $\mu_i^1+\mu_i^2+\lambda_i^1+\lambda_i^2,$ then given ${\tilde\omega},$ in $\tilde P_{\tilde\omega}$ distribution,
$$T_1=\sum_{k\le 0}\sum_{i=1}^{\overline U_k}\xi_{ki}.$$
Hence by Ward's equation, it follows from (\ref{enif}) that
\begin{eqnarray}\label{do}
 && \sum_{r=1}^2\tilde E_{\tilde\omega}(T_1|N_{T_1}=r)=\sum_{r=1}^2\sum_{k\le 0}\tilde E_{\tilde\omega}(\overline U_k|N_{T_1}=r)\tilde E_{\tilde\omega}(\xi_{k1})\no\\
&&= \sum_{k\le 0}\frac{1}{\mu_k^1+\mu_k^2+\lambda_k^1+\lambda_k^2}\no\\
&&\quad\quad\quad\times\Big(\frac{\alpha_{1,1}}{\alpha_{1,1}+\alpha_{1,2}},\frac{\alpha_{1,2}}{\alpha_{1,1}+\alpha_{1,2}},1,0,...,0\Big)\no\\
&&\quad\quad\quad\cdot(Q_0Q_{-1}\cdots
Q_{k+1}\mathbf v_1+Q_0Q_{-1}\cdots Q_{k}\mathbf v_2)\no\\
&&=:D({\tilde\omega})
\end{eqnarray}
and
\begin{equation}\label{pio}
  \begin{split}
  &\sum_{r=1}^2\sum_{k\le0}\tilde E_{\theta^{-k}{\tilde\omega}}\Big(\sum_{j=1}^{U_k}\xi_{kj}|N_{T_1}=r\Big)\\
&=\sum_{k= 0}^\infty\frac{1}{\mu_0^1+\mu_0^2+\lambda_0^1+\lambda_0^2}\\
&\quad\quad\quad\times\Big(\frac{\alpha_{k+1,1}}{\alpha_{k+1,1}+\alpha_{k+1,2}},\frac{\alpha_{k+1,2}}{\alpha_{k+1,1}+\alpha_{k+1,2}},1,0,...,0\Big)\\
&\quad\quad\quad\cdot(Q_kQ_{k-1}\cdots
Q_{1}\mathbf v_1+Q_kQ_{k-1}\cdots Q_{0}\mathbf v_2)\\
&=:\pi({\tilde\omega}).
\end{split}
\end{equation}
Substituting (\ref{do}) and (\ref{pio}) to (\ref{vp}),  we conclude that
\begin{equation*}
  \begin{split}
  v_{\tilde{\mathbb P}}&=\frac{\tilde\mbe \Big( \sum_{r=1}^2\sum_{k\le0}\tilde E_{\theta^{-k}{\tilde\omega}}\big(\sum_{j=1}^{U_k}\xi_{kj}|N_{T_1}=r\big)(2\lambda_0^2+\lambda_0^1-\mu_0^1-2\mu_0^2)\Big)}{\sum_{r=1}^R\tilde E(T_1|N_{T_1}=r)}\\
&=\frac{\tilde\mbe\z(\pi({\tilde\omega})(2\lambda_0^2+\lambda_0^1-\mu_0^1-2\mu_0^2)\y)}{\tilde\mbe(D({\tilde\omega}))}.
\end{split}
\end{equation*}
\qed

\section*{Appendix: On the existence of $\{N_t\}$}

Given $\tilde\omega,$ let $Q=(q_{ij})$ be a matrix with
$$q_{ij}=\left\{\begin{array}{ll}
                                       \lambda_i^r,&\text{ if } j=i+r,\ r=1,...,R; \\
                                       \mu_{i}^l,&\text{ if } j=i-l,\ l=1,...,L;\\
                                       -\big(\sum_{i=1}^L\mu_i^l+\sum_{r=1}^R\lambda_i^r), &\text{ if } j=i;\\
                                       0, &\text{ else.}
                                     \end{array}
\right.$$ Then  $Q$ is obviously a conservative Q-matrix. Note that under (C2), $Q$ is bounded from above. Hence the process $\{N_t\}$ exists (See for example Anderson \cite{ads}, Proposition 2.9, Chapter 2.). Next we give a condition which implies the existence of $\{N_t\}$ but is  weaker than (C2).
We have from classical argument that there exists at least one transition matrix $(\overline{p}_{\tilde\omega}(t,i,j))$ such that
\begin{equation}\label{uniq}
  \lim_{t\rightarrow0}\frac{\overline{p}_{{\tilde\omega}}(t,i,j)-\delta_{ij}}{t}=q_{ij},\ i,j\in\mathbb Z.
\end{equation}

Let $(p_{\tilde\omega}(h,i,j))$ be a standard transition matrix satisfying (\ref{uniq}). Let $\{N_t\}$ be a continuous time Markov chain with transition matrix  $(p_{\tilde\omega}(h,i,j)).$
Let  $\tau_0=0$ and
define $\tau_{n}:=\inf\{t>\tau_{n-1}:N_t\neq N_{\tau_{n-1}}\}$ recursively for $n\ge 1.$ Then $\tau_n,n\ge0$ are the consecutive discontinuities of the process $\{N_t\}.$
If \begin{equation*}\label{zz}\tilde\mbp\Big(\sum_{l=1}^L\mu_0^l+\sum_{r=1}^R\lambda_0^r>0\Big)=1,\end{equation*} then $\tilde P$-a.s., $\tau_n<\infty, n\ge0.$ Let $\chi_n=N_{\tau_n},$ for $n\ge0.$ The process $\{\chi_n\}_{n\ge0}$ is known as the embedded process of $\{N_t\}.$
\begin{proposition}
  Suppose  that $\tilde\mbp\big(\sum_{l=1}^L\mu_0^l+\sum_{r=1}^R\lambda_0^r>0\big)=1$ and
$$ \tilde{\mathbb P}\Big(\sum_{n=1}^\infty{\Big(\max_{1\le k\le R}\Big\{\sum_{r=1}^R\lambda_{nR-k}^r+\sum_{l=1}^L\mu_{nR-k}^l\Big\}\Big)^{-1}}=\infty\Big)=1,$$
 $$\tilde{\mathbb P}\Big(\sum_{n=-\infty}^0{\Big(\max_{1\le k\le L}\Big\{\sum_{r=1}^R\lambda_{nL-k}^r+\sum_{l=1}^L\mu_{nL-k}^l\Big\}\Big)^{-1}}=\infty\Big)=1.$$ Then for $\tilde\mbp$-a.a. ${\tilde\omega},$ there is a unique transition matrix $(p_{\tilde\omega}(h,i,j))$ which satisfies (\ref{uniq}).
\end{proposition}
\proof Considering the above defined $\{N_t\},$  since $\tilde\mbp\big(\sum_{l=1}^L\mu_0^l+\sum_{r=1}^R\lambda_0^r>0\big)=1,$ we have $\tilde P$-a.s., $\tau_n<\infty, n\ge 1.$ By the classical argument of the uniqueness of the $Q$-process, if $$\tilde P(\lim_{n\rto}\tau_n=\infty)=1,$$ then the minimal solution $p_{{\tilde\omega}}(t,i,j)$ is the unique $Q$-transition matrix. Let $q_{i}=-q_{ii},\ i\in\mbz.$
If \begin{equation}\label{tin}
  \tilde P\Big( \sum_{n=0}^\infty q^{-1}_{\chi_n}=\infty\Big)=1,
\end{equation} then we have (see Chung \cite{chung}, Theorem 1 in II.19) that $\tilde P(\lim_{n\rto}\tau_n=\infty)=1.$
Next we show (\ref{tin}). In fact, if the process $\{\chi_n\}_{n\ge 0}$ is recurrent or transient to the right, it must visit at least one state of each of the sets $B_n:=\{nR-k\}_{k=1}^R,n=1,2,....$ Then $\tilde P$-a.s.,
$$\sum_{n=0}^\infty q^{-1}_{\chi_n}\ge \sum_{n=1}^\infty{\Big(\max_{1\le k\le R}\Big\{\sum_{r=1}^R\lambda_{nR-k}^r+\sum_{l=1}^L\mu_{nR-k}^l\Big\}\Big)^{-1}}=\infty.$$
Else if the process $\{\chi_n\}_{n\ge 0}$ is transient to the left, it must visit at least one state of each of the sets $A_n:=\{nL-k\}_{k=1}^L,n=0,-1,-2,....$  It follows that $ P$-a.s., $$\sum_{n=0}^\infty q^{-1}_{\chi_n}\ge \sum_{n=-\infty}^0{\Big(\max_{1\le k\le L}\Big\{\sum_{r=1}^R\lambda_{nL-k}^r+\sum_{l=1}^L\mu_{nL-k}^l\Big\}\Big)^{-1}}=\infty.$$ Consequently (\ref{tin}) follows.\qed

\noindent{\large{\bf \large Acknowledgements:}} The author would like to thank Professor Wenming Hong for his useful comments on the paper.

 % \begin{center}
%{\section*{Acknowledgements}}
%\end{center}

\end{document}